\documentclass[runningheads,a4paper]{llncs}

\usepackage{amsmath,amssymb}
\usepackage{mathtools}

\usepackage{graphicx}
\usepackage{color}
\usepackage{bm}
\usepackage{subcaption}

\usepackage{siunitx}

\newcommand{\R}{\mathbb{R}}

\usepackage{url}
\usepackage[bookmarks=false,hyperfootnotes=false]{hyperref}
\hypersetup{
			colorlinks=true,
			linkcolor=red,
			anchorcolor=black,
			citecolor=green,
			urlcolor=blue
			}
\begin{document}

\mainmatter  

\title{Digital Cultural Heritage imaging \\ via osmosis filtering} 
\titlerunning{Digital Cultural Heritage imaging via osmosis filtering}

\author{Simone Parisotto\inst{1} \and Luca Calatroni\inst{2} \and Claudia Daffara\inst{3}}
\authorrunning{S. Parisotto, L. Calatroni \and C. Daffara}

\institute{
CCA, Wilberforce Road, CB3 0WA, University of Cambridge, UK\\
\email{sp751@cam.ac.uk}
\and
CMAP, Ecole Polytechnique 91128 Palaiseau Cedex, France\\
\email{luca.calatroni@polytechnique.edu},
\and
University of Verona, Strada Le Grazie 15, 37134 Verona, Italy\\
\email{claudia.daffara@univr.it}
}

%
%

\toctitle{Lecture Notes in Computer Science}
\tocauthor{Authors' Instructions}
\maketitle

\begin{abstract} 
In Cultural Heritage (CH) imaging, data acquired within different spectral regions are often used to inspect surface and sub-surface features. 
Due to the experimental setup, these images may suffer from intensity inhomogeneities, which may prevent conservators from distinguishing the physical properties of the object under restoration. 
Furthermore, in multi-modal imaging, the transfer of information between one modality to another is often used to integrate image contents.

In this paper, we apply the image osmosis model proposed in \cite{weickert,vogel,Calatroni2017} 
to solve correct these problems arising when diagnostic CH imaging techniques based on reflectance, emission and fluorescence mode in the optical and thermal range are used.
For an efficient computation, we use stable operator splitting techniques to solve the discretised model. We test our methods on real artwork datasets: the thermal measurements of the mural painting ``Monocromo'' by Leonardo Da Vinci, the UV-VIS-IR imaging of an ancient Russian icon and the Archimedes Palimpsest dataset. 

\end{abstract}

\section{Introduction}
Non-destructive imaging techniques are a well established tool in the diagnostics of artworks. Different spectral bands may unveil hidden information about the surface and the sub-surface layers. 
Most notable methods within the optical range include UV-based imaging in an emission mode, which exploits the fluorescence of materials to map pigments and past restorations, and IR-based imaging in reflectance mode (reflectography), which allows the detection of underlying features (e.g.\ the artist preparatory drawing and repaintings) thanks to the optical transparency of pictorial matter in the near-IR (from \SIrange{0.8}{2.5}{\micro\meter}). 
Recently, it has been shown that reflectography techniques can be performed also in the mid-IR thermal region (from \SIrange{3}{5}{\micro\meter}) and that such method, called Thermal-Quasi Reflectography (TQR), can complement  optical techniques in the analysis of frescoes \cite{Daffara2012,Daffara2017}.  

Artwork analysis demands the above procedures to have sub-millimetric spatial resolution and full-field capability, 
which often (large paintings, small sensor array) leads to the practice of mosaicking images acquired in multiple acquisition sessions, each one being the result of its specific setup. 
Irradiance inhomogeneities, for instance, can be induced by a number of experimental factors related to the kind of acquisition modality, source-camera positioning as well as environmental light. The intensity values in the final mosaic are thus affected by errors and must be post-processed by suitable radiometric procedures in order to obtain the correct reflectance or fluorescence map. 

In this work, we discuss the use of the linear image osmosis model proposed in \cite{weickert,vogel,Calatroni2017} as an unsupervised mathematical tool to correct the inter-frame light differences and to transfer image information in multi-modal imaging \cite{CalatroniUnveiling}.
Our numerical implementation relies on the use of mathematical operator splitting techniques considered, e.g., in \cite{Barash2001,hundsbook,CalatroniADI,Calatroni2017}, which allows significant computational speed-ups for large collages. 
Numerical methods are validated on  real artwork datasets acquired via different imaging techniques:
a TQR mosaic of the ``Monocromo" by Leonardo Da Vinci (Castello Sforzesco, Milan, Italy) \cite{Daffara2015}, a stack of UV-VIS-IR images of an ancient Russian icon and multi-modal data of the Archimedes palimpsest.

\paragraph{Structure of the Paper.} 
In Section \ref{sec: challenges} we present some of the challenges frequently encountered in CH imaging which we are going to address, with a brief review of the image modalities used to acquire the datasets used for our tests. 
In Section \ref{sec: linear_osmosis} the mathematical model of image osmosis  \cite{vogel,weickert} is revised and for its numerical realisation operator splitting techniques are briefly described. 
In Section \ref{sec: numerical} we report the results obtained using the osmosis model to solve the challenges described.

\section{Intensity balance in CH imaging modalities}\label{sec: challenges} 

\paragraph{Optical Ultraviolet (UV)- Visible (VIS) - Infrared (IR) imaging.}
CH imaging in the optical region from \SIrange{0.35}{1.1}{\micro\meter} is performed with 
CCD or CMOS sensors with state-of-the-art arrays up to 100Mp. Multimodal UV-VIS-IR images are acquired in fix geometry, changing lamps and camera setting, and combined in a stack for their subsequent analysis.
Problems of large-area uniform irradiance are typical due to the geometry and the size of the available sources, e.g. wood tubes, bulbs, flash, led arrays. 

IR reflectography records the near-IR radiation 
reflected by the object after scattering through the painting layers: image intensity values are here related to the presence of an absorbing/reflective medium. Therefore, the intensity map allows the description of sub-surface features by contrast. 
In an ideally correct UV fluorescence image, the intensities represent the signal emitted in the VIS range in response to UV-light stimulus and are related to the kind and quantity of the materials in the outer surface. 

In multimodal imaging the raw signal is collected using linear detectors for all the different modalities considered, but the acquisition procedure is different (reflectance, thermal quasi-reflectance, luminescence) as well as the employed radiation (optical and thermal). 

In CH diagnostics, the fusion of multimodal optical images in pseudo-colour schemes is typically used to enhance differences in the spectral reflectance. A notable example is the IR False Colour technique (IR-R-G), which allows discriminating pigments similar in the VIS but with different IR behaviour. 

\paragraph{Thermal-Quasi Reflectography (TQR).}\label{sec: tqr}
TQR is a recent imaging technique for analysing frescoes \cite{Daffara2012} 
based on the acquisition of 
the thermal mid-IR (from \SIrange{3}{5}{\micro\meter}) signal in reflectance mode. 
The 
rational behind
quasi-reflectography is that an object at room temperature 
emits only $1.1\%$ of its thermal energy in the mid-IR: thus, 
a signal dominated by the reflected thermal radiation can be collected
by sending a non-heating mid-IR stimulus.

In TQR imaging, the intensities in the final image represent the TQR reflectance $r$, defined as the ratio of the incident and reflected radiant flux in the mid-IR, and obtained point-wise from the acquired radiometric image $u$ using an in-scene calibration target of known reflectance $r_\mathrm {ref}$, i.e. $
u_{ij} = r_{ij} u_{\mathrm{ref}} r_{\mathrm{ref}}^{-1},
$
where $u_{\mathrm{ref}}$ is the averaged radiometric image of the target. 

Thermal camera are equipped with small array sensors (1Mp). The uniform irradiance of the object is often hard to achieve for thermal sources due to reflections from the environment. 
Furthermore, an increase of the object surface temperature may occur during the time of a multiple-frames measuring session, causing a variation of the emitted radiance and a shift of the spectral distribution towards the mid-IR \cite{Daffara2017}, which results in inter-frame dis-homogeneity.  
As a consequence, the mosaic of the TQR frames is often affected by piecewise non-uniform illumination, even after the calibration with the reference target. 

\paragraph{Quantitative material mapping.} In all the imaging modalities described above, it is therefore quite common observing inhomogeneities in the acquired intensity image depending on the experimental setup and measurement errors. A careful processing of the resulting intensity image is therefore desirable to quantify and map the composition of the material under consideration. To solve this problem we consider in the following a mathematical drift-diffusion PDE model.

\section{Image Osmosis Filtering} \label{sec: linear_osmosis}
For a regular image domain $\Omega \subset\R^2$ with boundary $\partial \Omega$ and given a vector field $\bm{d}:\Omega\to\R^2$, the linear image osmosis model considered in \cite{vogel,weickert,Calatroni2017} is a drift-diffusion PDE which computes for every $t\in (0,T],~ T>0$, a family $\left\{u(x,t)\right\}_{t>0}$ of regularised images of a positive initial gray-scale image $f:\Omega\to\R^+$ by solving:
\begin{equation} \label{eq:osmosis}
\begin{cases}
\partial_t u = \Delta u - \text{div}(\bm{d}u) & \text{ on } \Omega \times (0,T] \\
u(x,0)=f(x) & \text{ on } \Omega \\
\langle \nabla u - \bm{d}u, \bm{n} \rangle = 0 & \text{ on } \partial \Omega \times (0, T],
\end{cases}
\end{equation}
where $\langle \cdot, \cdot \rangle$ is the Euclidean scalar product and $\bm{n}$ the outer normal vector on $\partial\Omega$. Any solution of \eqref{eq:osmosis} preserves the average gray value (AVG) of $f$ and it is non-negative at any time $t>0$, see \cite[Proposition 1]{weickert}. Moreover, by setting $\bm{d}:=\bm{\nabla} \ln v$ for a given reference image $v>0$, the steady state of \eqref{eq:osmosis} is a rescaled version of $v$, i.e. $w(x)=\frac{\mu_f}{\mu_v}~v(x)$, where $\mu_f$ and $\mu_v$ are the average of $f$ and $v$ over $\Omega$, respectively. To avoid the dependence on the reference image $v$ (not available in practice), the vector field $\bm{d}$ can be modified 
to balance intensity information in the image or to integrate on one image the information contained in another image (seamless image cloning).

\subsection{Numerical realisation} 
In \cite{vogel}, a fully discrete theory for the model \eqref{eq:osmosis} is studied. For a given  $\bm{f}\in\R^N_+$, a spatial finite-difference discretisation of the differential operators in \eqref{eq:osmosis} reads
\begin{equation}   \label{eq:semi_discrete}
\bm{u}(0) = \bm{f}, \qquad \bm{u}'(t) =\bm{A} \bm{u}(t),\qquad t>0.
\end{equation}

For the time discretisation of \eqref{eq:semi_discrete},
in \cite{vogel,weickert} the authors consider standard forward and backward Euler schemes and prove conservation and convergence properties analogous to the continuous model \eqref{eq:osmosis}, \cite[Proposition 1]{vogel}. In particular, unconditional stability is observed for the fully implicit case.  Computationally, iterative solvers (e.g., BiCGStab) are used to solve the resulting penta-diagonal and non-symmetric linear systems. 

\paragraph{Operator splitting schemes.} 
The computation of the numerical solution of \eqref{eq:semi_discrete} may be extremely costly for large images. To overcome this issue, we consider in the following a method similar to the one proposed in \cite{Calatroni2017} and based on a dimensional splitting idea.

Given a domain $\Omega\in\R^s$ for $s\geq 1$, \emph{dimensional splitting} methods are based on the decomposition of the operator $\bm{A}$ of the initial boundary value problem \eqref{eq:semi_discrete} into the sum:
\begin{equation} \label{decomp}
\bm{A}=\bm{A_0}+\bm{A_1}+\ldots+\bm{A_s}
\end{equation}
where the terms $\bm{A_j}$, $j=1,\ldots, s,$ encode the linear action of $\bm{A}$ 
along the space direction $j$ while $\bm{A_0}$ may contain mixed-derivative, non-stiff and non-linear terms.
 
Note that when $\bm{A_0}=0$ additive (AOS) and multiplicative (MOS) operator splitting techniques can be used to solve the problem in a fully implicit way, thus ensuring unconditional stability \cite{hundsbook}.

Using standard finite difference space discretisation, such methods requires the numerical inversion of tridiagonal matrices, which can be rendered efficiently via standard matrix factorization techniques.

In our case $s=2$, $\bm{A}_0=0$ and the tridiagonal matrices $\bm{A_1}$, $\bm{A_2}$ read 
\begin{align}  \label{eq:discr_operators_splitted}
\bm{A_1}\bm{u}&  = u_{i,j} \bigg( -\frac{2}{h^2}+\frac{d_{1,i-\frac{1}{2},j}}{h}-\frac{d_{1,i+\frac{1}{2}}}{h}+\frac{d_{2,i,j-\frac{1}{2}}}{h}-\frac{d_{2,i,j+\frac{1}{2}}}{h} \bigg) \notag \\
& +  u_{i+1,j} \bigg( \frac{1}{h^2}-\frac{d_{1,i+\frac{1}{2},j}}{2h} \bigg)+ u_{i-1,j} \bigg( \frac{1}{h^2}+\frac{d_{1,i-\frac{1}{2},j}}{2h} \bigg);\\
\bm{A_2}\bm{u}&  = u_{i,j} \bigg( -\frac{2}{h^2}+\frac{d_{1,i-\frac{1}{2},j}}{h}-\frac{d_{1,i+\frac{1}{2}}}{h}+\frac{d_{2,i,j-\frac{1}{2}}}{h}-\frac{d_{2,i,j+\frac{1}{2}}}{h} \bigg) \notag \\
& + u_{i,j+1} \bigg( \dfrac{1}{h^2}-\dfrac{d_{2,i,j+\frac{1}{2}}}{2h} \bigg)+
u_{i,j-1} \bigg( \dfrac{1}{h^2}+\dfrac{d_{2,i,j-\frac{1}{2}}}{2h} \bigg). \notag
\end{align}
We make use of the Additive Operator Splitting (AOS) numerical scheme \cite{hundsbook} where for every iteration $k\geq 1$ and any time-step $\tau>0$ the update reads
\begin{align}  
\mathbf{u}^{k+1} &= \frac{1}{2}~\sum_{n=1}^2~\left(\mathbf{I} -2\tau \mathbf{A}_n\right)^{-1}\mathbf{u}^{k}. 
\tag{AOS} 
\label{eq:AOS}
\end{align}
The AOS method is unconditionally stable and first-order accurate in time.  In \cite{Parisotto2018} the scale-space properties of such splitting scheme applied to the osmosis problem are shown. 

\section{Applications}\label{sec: numerical}
In the following we adapt the osmosis model \eqref{eq:osmosis} for the CH imaging problems described in Section \ref{sec: challenges}. 

\paragraph{Datasets.} 
The TQR mosaic in Figure \ref{fig: TQR mosaic} ($4717\times 7066$ pixels) is made by 33 overlapped TQR tiles shot by FLIR~X6540sc camera used as a radiometer on the ``Monocromo'' wall fresco (L.\ Da Vinci, Castello Sforzesco, Milan, Italy) in Figure \ref{fig: visible orthophoto}, see \cite{Daffara2015,Daffara2017} for more details on the workflow.
The UV and IR data acquired on a Russian icon in Figure \ref{fig: visible orthophoto} compose a mosaic of 5 overlapped tiles ($3671\times 4551$ pixels) shot by a Nikon D810A camera.
The Archimedes Palimpsest in Figure \ref{fig: arch input} ($1200\times 1600$ pixels) reveals a 10th-century copy of works by Archimedes of Syracuse in Figure \ref{fig: arch written}, a pseudo-colour UV-RGB separation.

\subsection{Light and colour correction}\label{sec: light correction}
\begin{figure}
\centering
\begin{subfigure}{0.3\textwidth}
\includegraphics[width=1\textwidth,trim=0 0 0 0cm,clip=true]{\detokenize{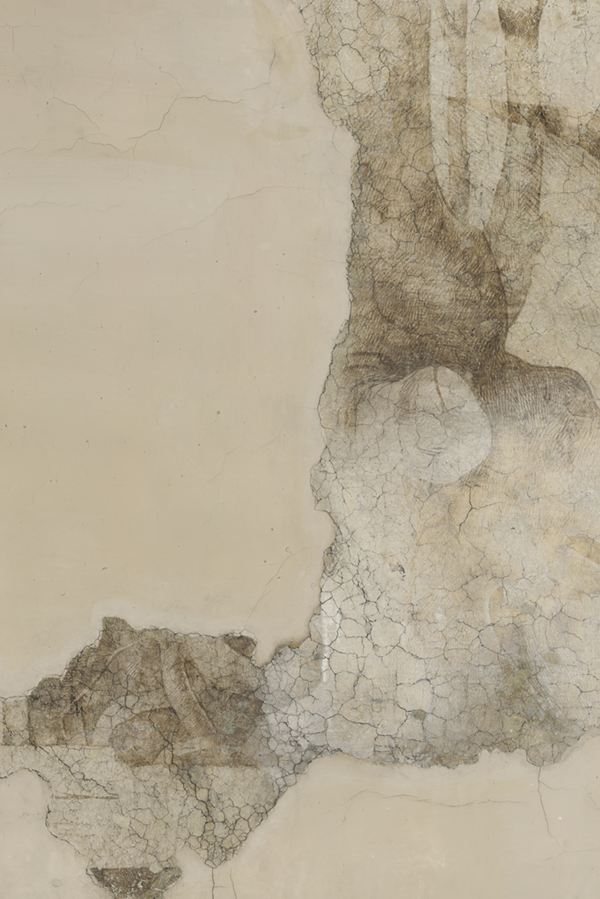}}
\caption{Visible Orthophoto.}
\label{fig: visible orthophoto}
\end{subfigure}
\hfill
\begin{subfigure}{0.3\textwidth}
\includegraphics[width=1\textwidth,trim=0 0 0 0cm,clip=true]{\detokenize{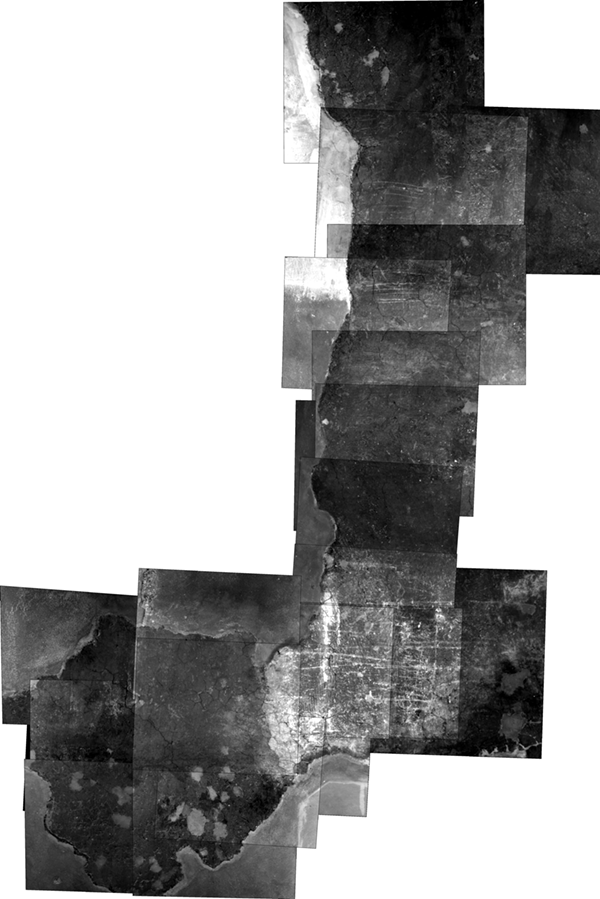}}
\caption{TQR mosaic.}
\label{fig: TQR mosaic}
\end{subfigure}
\hfill
\begin{subfigure}{0.3\textwidth}
\includegraphics[width=1\textwidth,trim=0 0 0 0cm,clip=true]{\detokenize{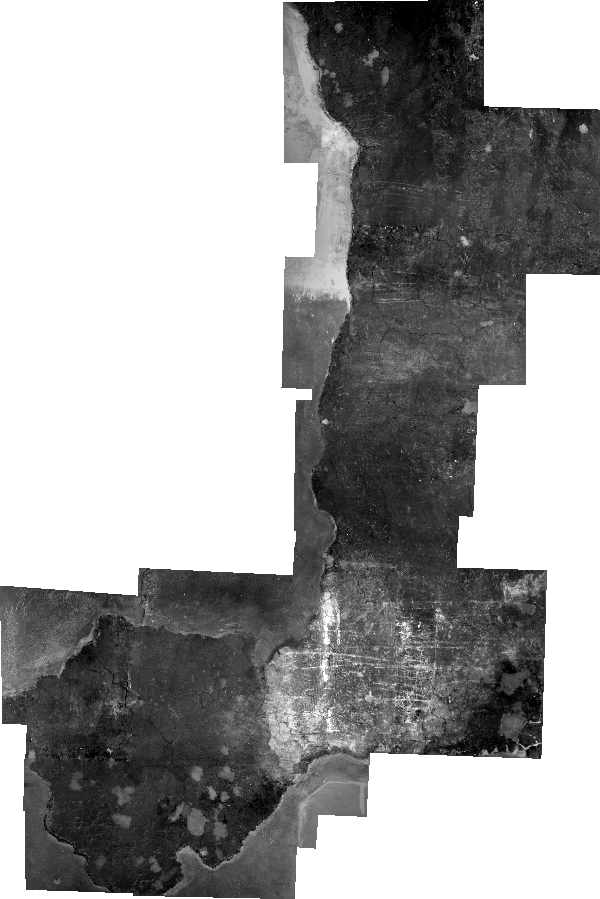}}
\caption{Osmosis result (\SI{629}{\second}).}
\label{fig: TQR AOS}
\end{subfigure}
\caption{Light-balance on TQR data: $T=1\mathrm{e}5$, $\tau=1\mathrm{e}3$.}
\label{fig: TQR BEFORE AFTER}
\end{figure}

We start testing the osmosis model on the light balance problems described in Section \ref{sec: challenges} for TQR and IR imaging. Let $\Omega_{b}\subset\Omega$ be the boundary between the collaged images, for this application the vector field $\bm{d}$ is defined as: $\bm{d} = \chi_{\Omega\setminus\Omega_b} \nabla\ln f$. In Figure \ref{fig: TQR BEFORE AFTER} we consider a TQR mosaic of images of the ``Monocromo" fresco showing light inhomogeneities. We compute the solution by \eqref{eq:AOS} scheme with nested \texttt{LU} factorisation of the tridiagonal matrices. 
We fix the parameters as $T=1\mathrm{e}5$, with time-step $\tau=1000$. The elapsed time is reported in the caption.

\begin{figure}
\centering
\begin{subfigure}{0.225\textwidth}\centering
\includegraphics[width=0.95\textwidth]{\detokenize{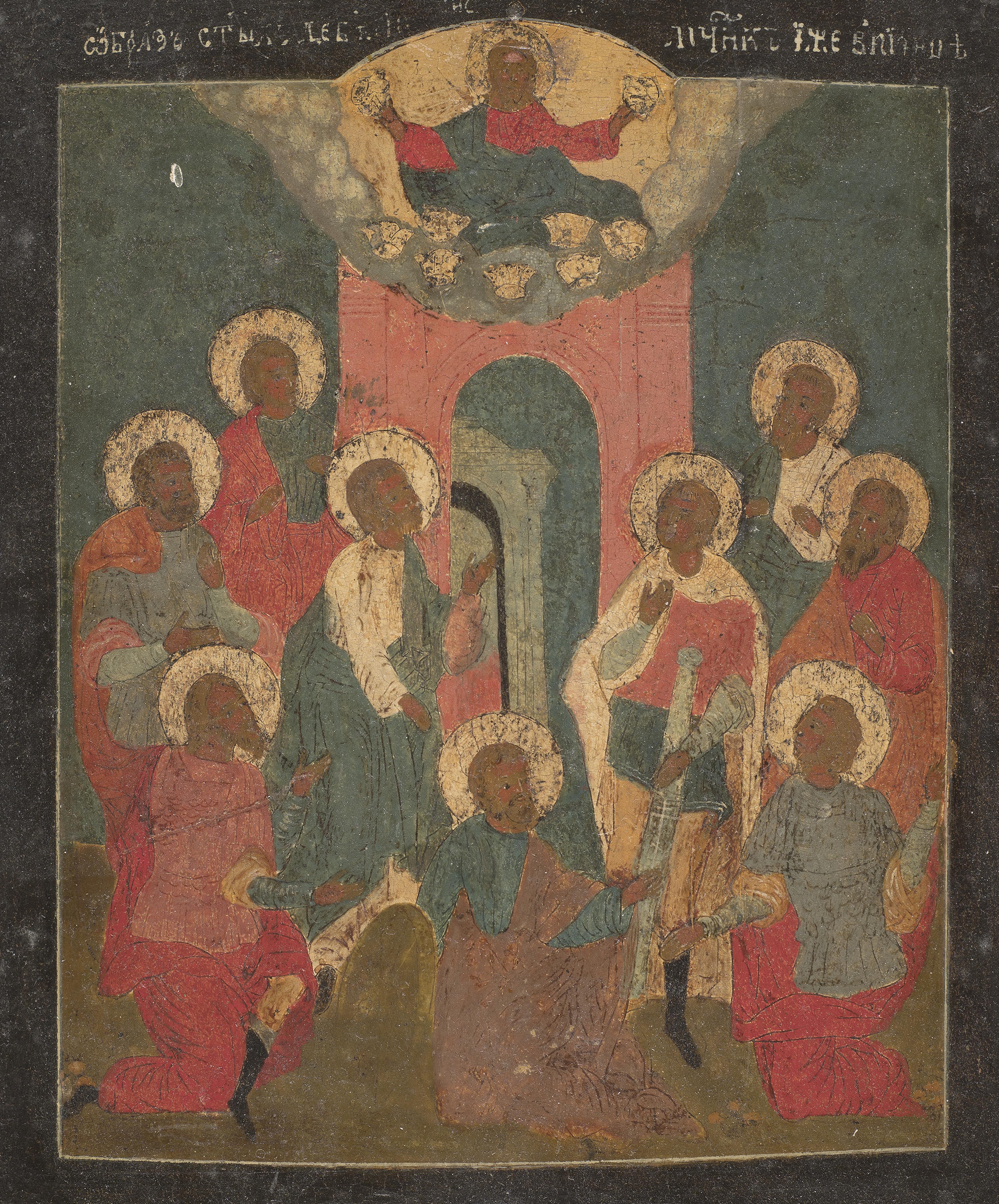}}
\caption{Full Visible}
\label{fig: russian icon1}
\end{subfigure}
\hfill
\begin{subfigure}{0.378\textwidth}
\centering
\includegraphics[width=1\textwidth,trim=0 0cm 0 0cm,clip=true]{\detokenize{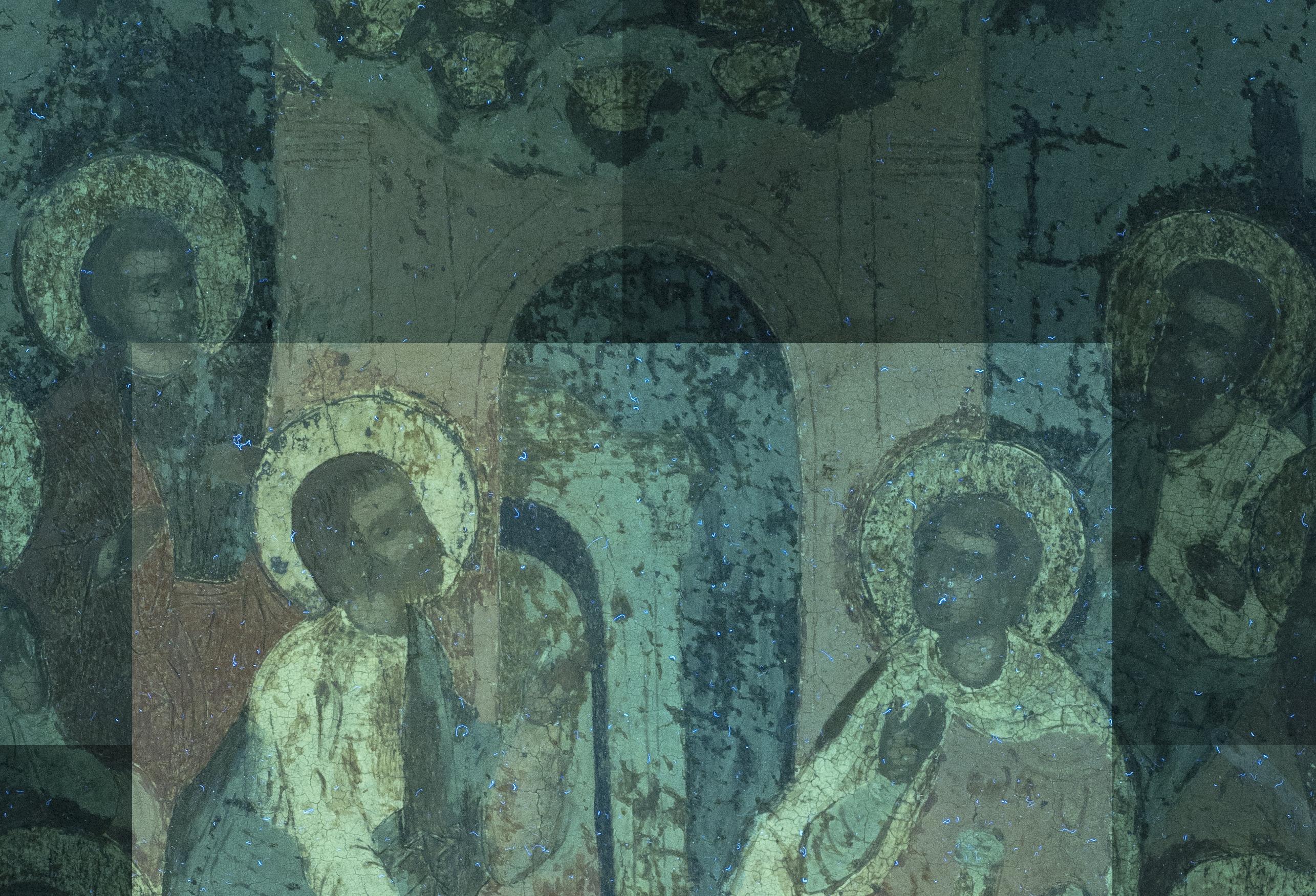}}
\caption{UV Fluorescence}
\label{fig: uvf pre}
\end{subfigure}
\hfill
\begin{subfigure}{0.378\textwidth}
\centering
\includegraphics[width=1\textwidth,trim=0 0cm 0 0cm,clip=true]{\detokenize{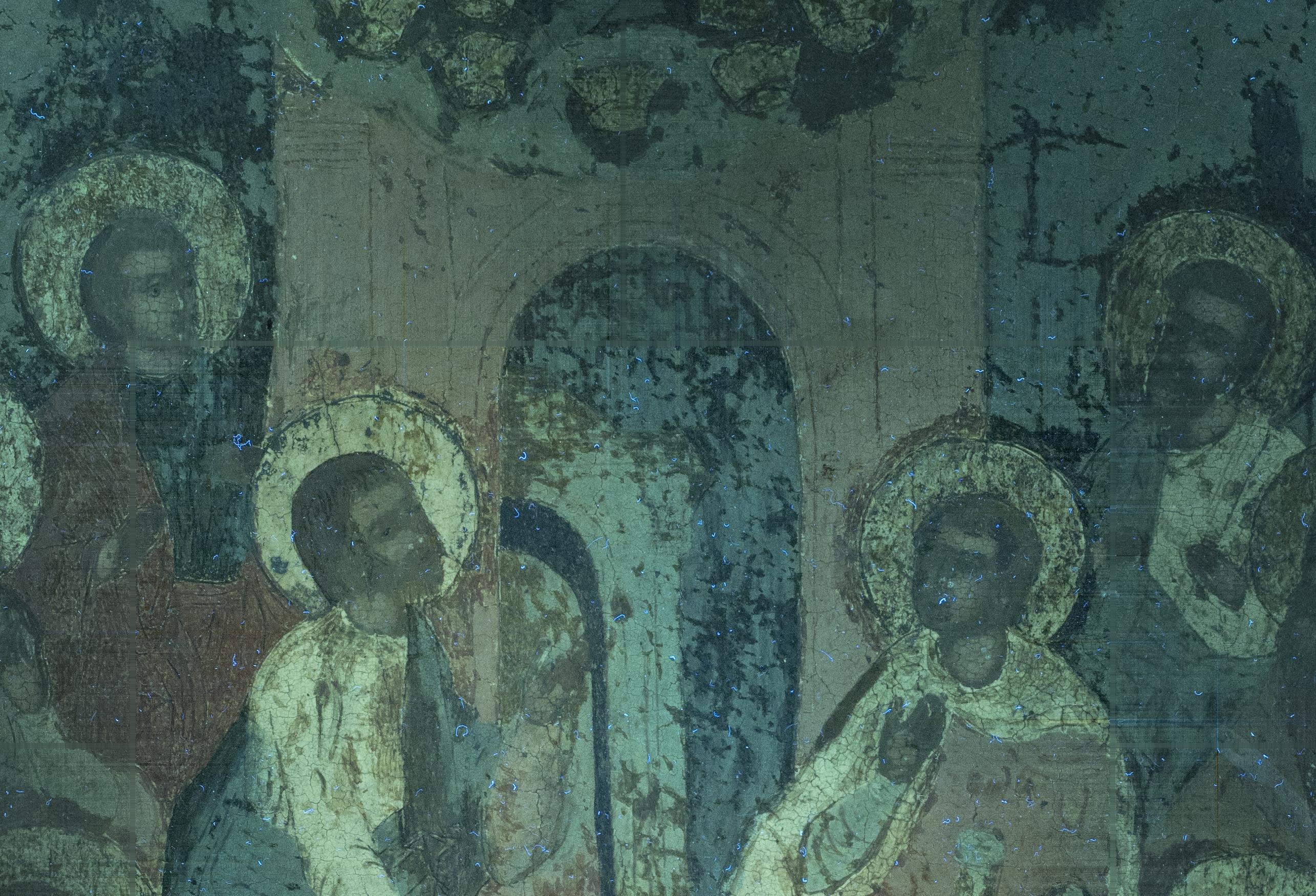}}
\caption{Osmosis result (1693 s.)}
\label{fig: uvf post}
\end{subfigure}
\caption{Light-balance on UV Fluorescence data (detail): $T=1\mathrm{e}5$, $\tau=1\mathrm{e}3$.}
\label{fig: uv fluo}
\end{figure}

A similar approach is applied in Figure \ref{fig: uv fluo} for UV fluorescence imaging data sampled from the Russian icon in Figure \ref{fig: russian icon1}. Once again the inhomogeneities due to the collage of multiple samples (inter-frame error) are minimized. 

\subsection{False colour adjustment}
\begin{figure}
\centering
\begin{subfigure}{0.49\textwidth}
\centering
\includegraphics[width=1\textwidth,trim=0 3cm 0 15cm,clip=true]{\detokenize{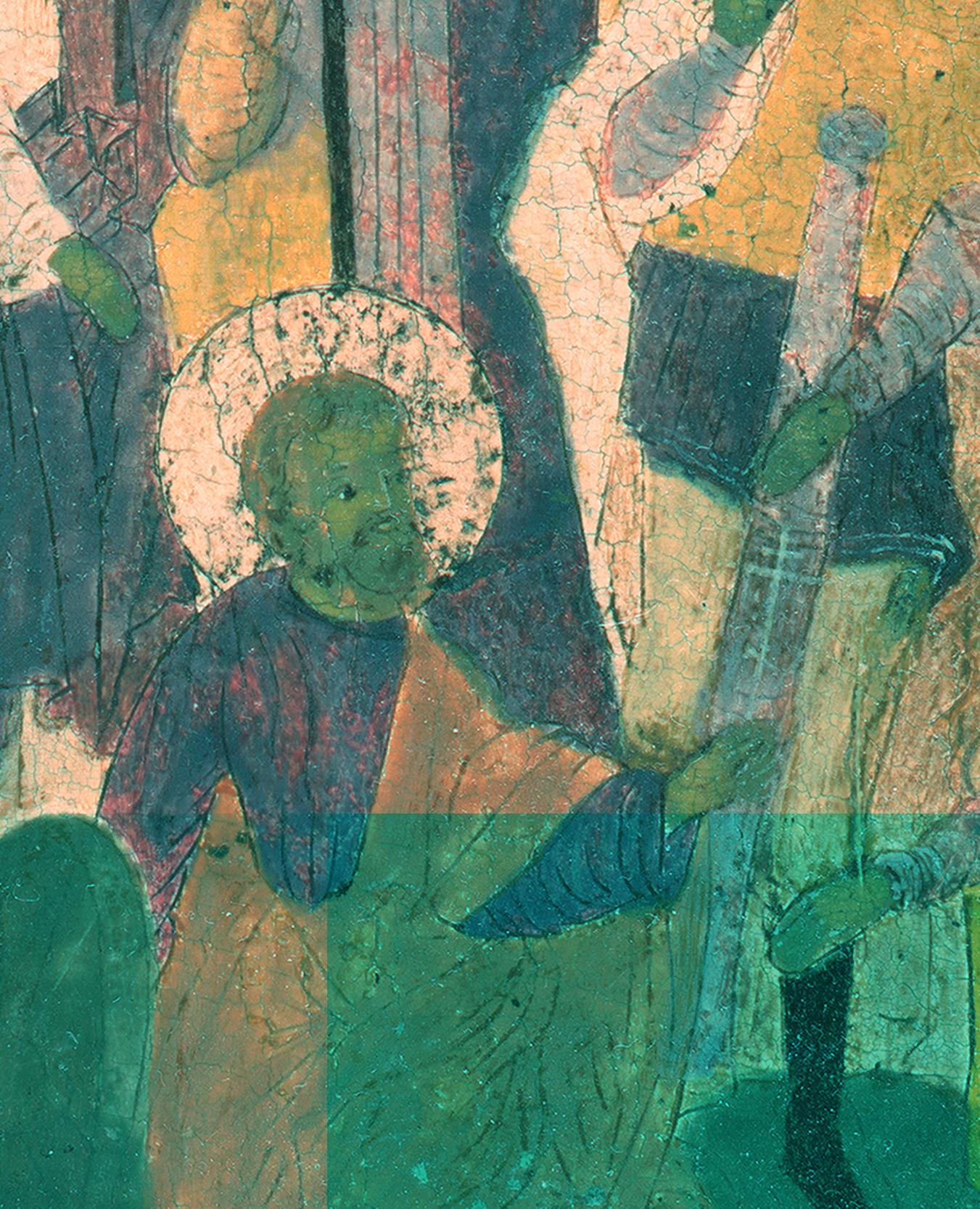}}
\caption{False-colour (input IR)-R-G}
\label{fig: falsecolor vis2inp}
\end{subfigure}
\hfill
\begin{subfigure}{0.49\textwidth}
\centering
\includegraphics[width=1\textwidth,trim=0 3cm 0 15cm,clip=true]{\detokenize{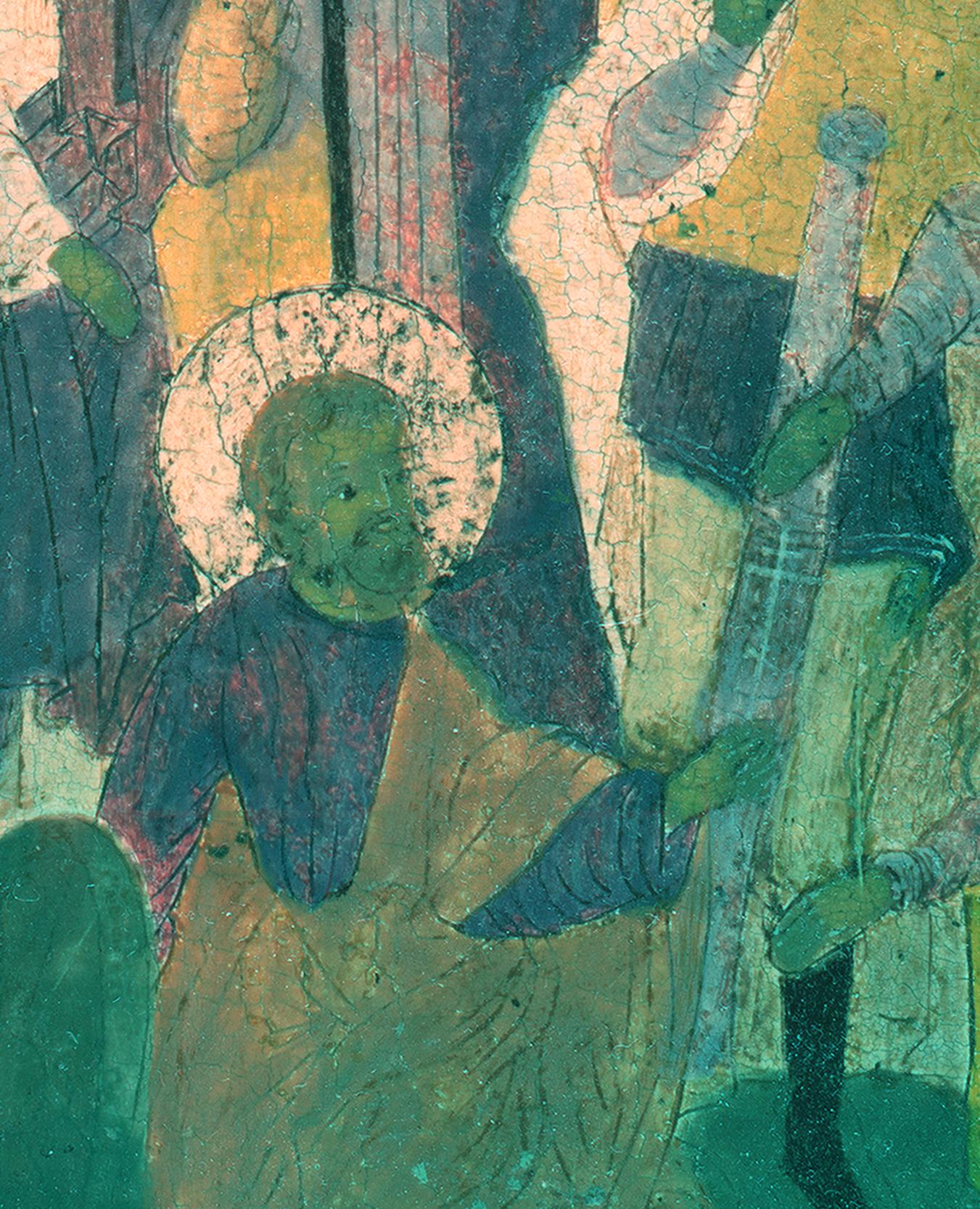}}
\caption{False-colour IR-R-G (721 s.)}
\label{fig: falsecolor vis2out}
\end{subfigure}
\caption{Light-balance for (zoomed) false-colour images: $T=1\mathrm{e}5$, $\tau=1\mathrm{e}3$.}
\label{fig: falsecolor}
\end{figure}

In Figure \ref{fig: falsecolor}, we observe how the measurement errors affect the false colour rendering of the Russian icon in Figure \ref{fig: russian icon1} and lead to incorrect material mapping. In particular, note in Figure \ref{fig: falsecolor vis2inp} the false colour IR-R-G image where the IR channel is affected by measurement errors. The use of the osmosis model for IR light balance makes the different regions  comparable, see Figure \ref{fig: falsecolor vis2out}.

\subsection{Multi-modal data integration}

\begin{figure}
\centering
\begin{subfigure}{0.31\textwidth}
\centering
\includegraphics[width=1\textwidth,trim=0 4.4cm 0 1.2cm,clip=true]{\detokenize{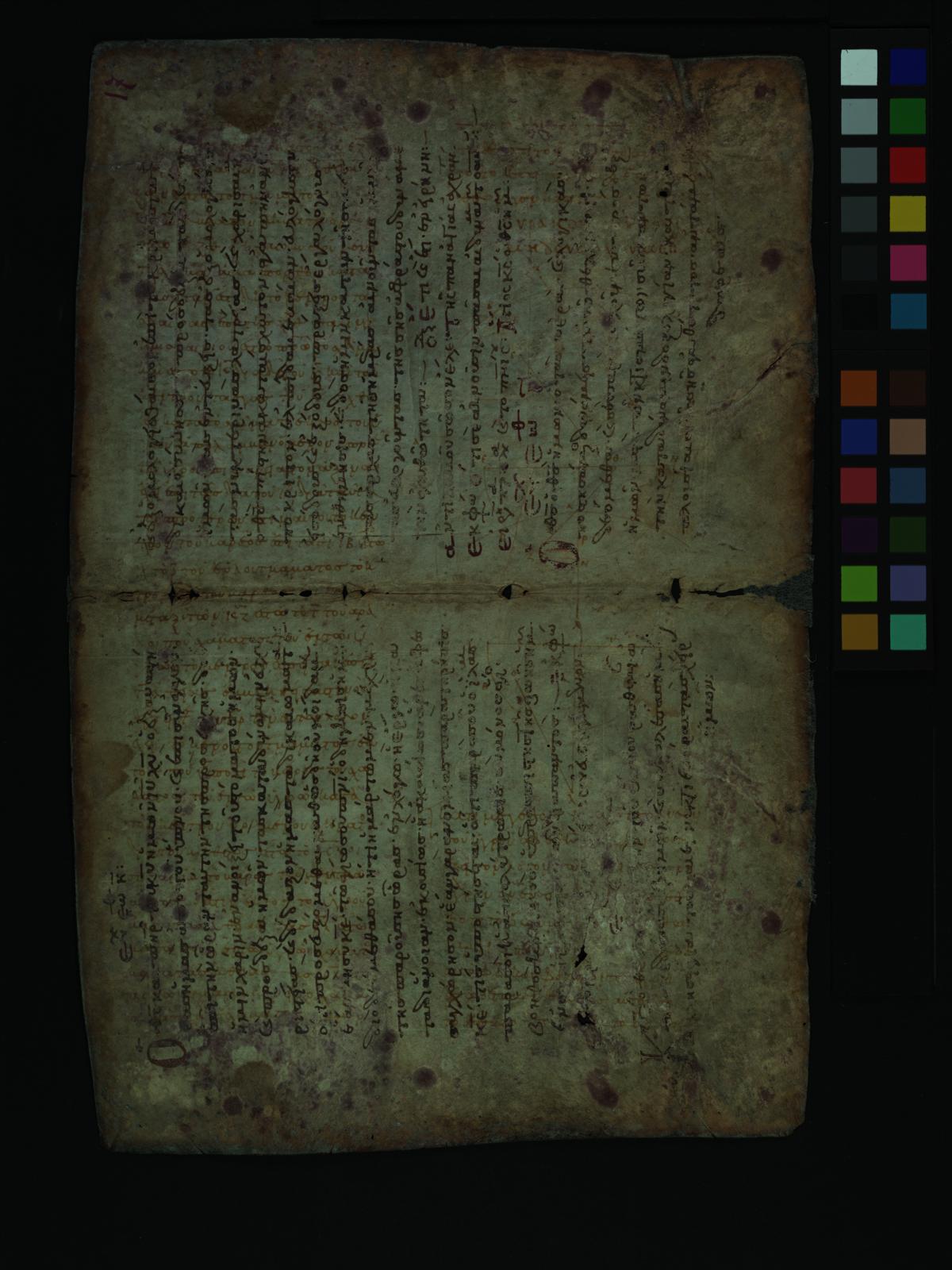}}
\end{subfigure}
\hfill
\begin{subfigure}{0.31\textwidth}
\centering
\includegraphics[width=1\textwidth,trim=0 4.4cm 0 1.2cm,clip=true]{\detokenize{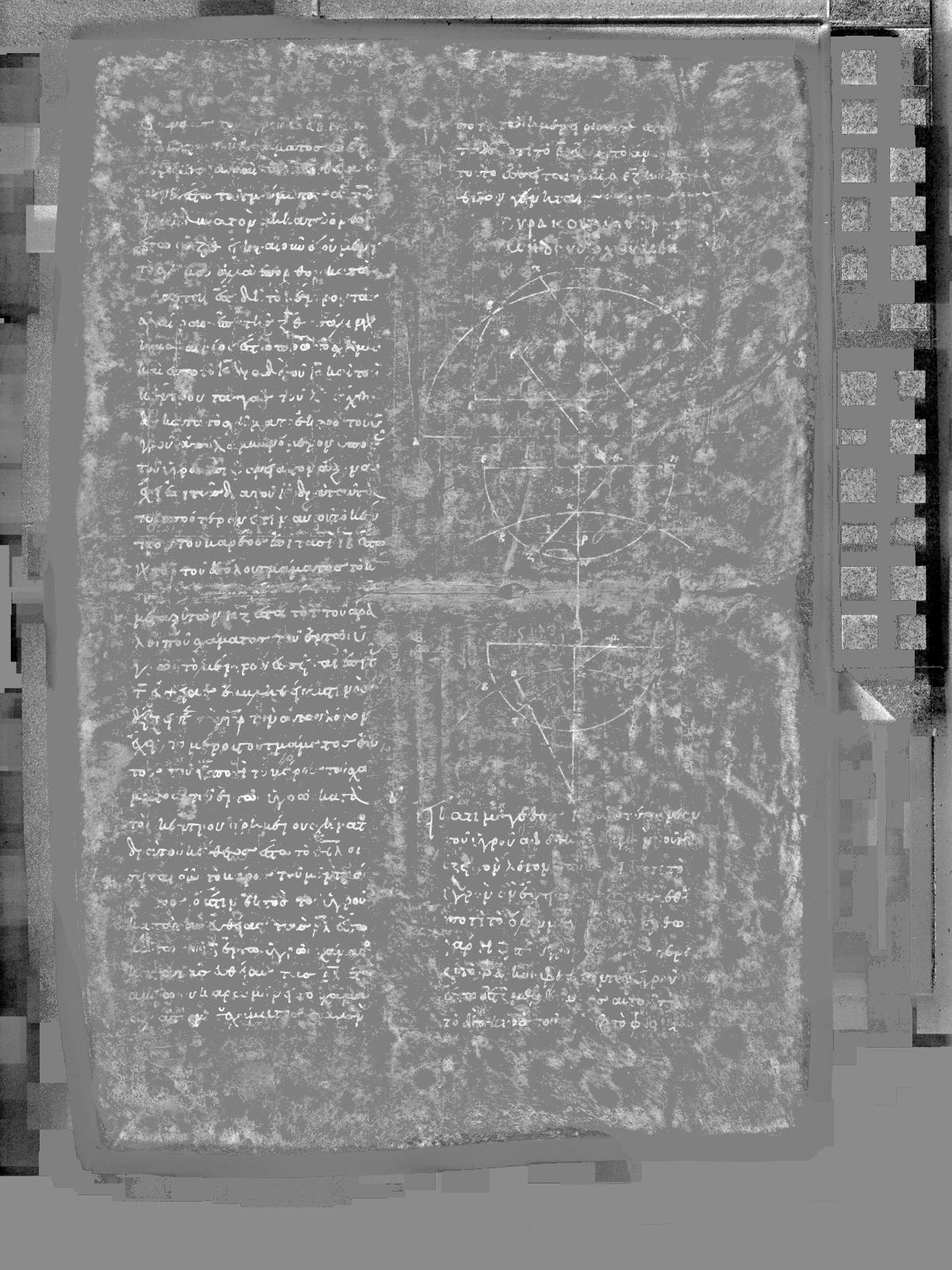}}
\end{subfigure}
\hfill
\begin{subfigure}{0.31\textwidth}
\centering
\includegraphics[width=1\textwidth,trim=0 4.4cm 0 1.2cm,clip=true]{\detokenize{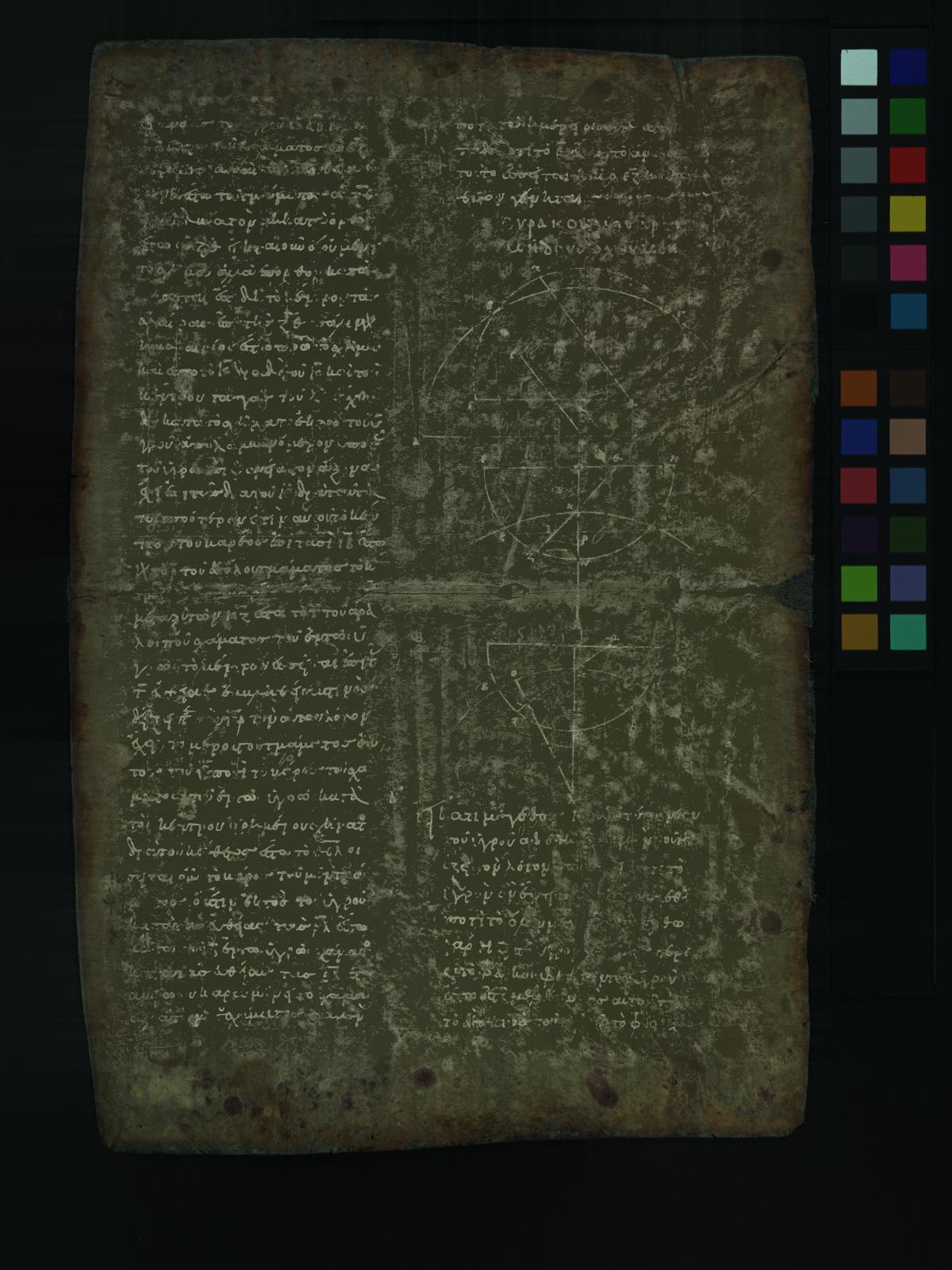}}
\end{subfigure}
\\\vspace{0.2em}
\begin{subfigure}{0.31\textwidth}
\centering
\includegraphics[width=1\textwidth,trim=0 2.3cm 0 5cm,clip=true]{\detokenize{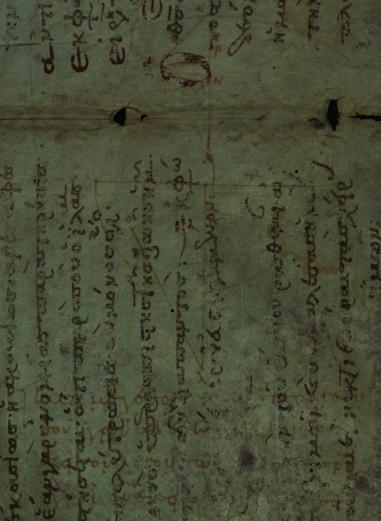}}
\caption{Input}
\label{fig: arch input}
\end{subfigure}
\hfill
\begin{subfigure}{0.31\textwidth}
\centering
\includegraphics[width=1\textwidth,trim=0 2.3cm 0 5cm,clip=true]{\detokenize{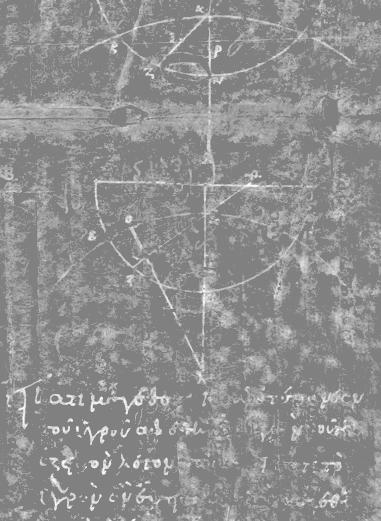}}
\caption{Pseudo-colour $v$}
\label{fig: arch written}
\end{subfigure}
\hfill
\begin{subfigure}{0.31\textwidth}
\centering
\includegraphics[width=1\textwidth,trim=0 2.3cm 0 5cm,clip=true]{\detokenize{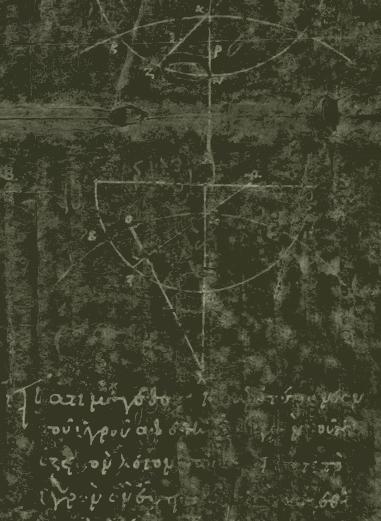}}
\caption{Osmosis Result (137 s.)}
\end{subfigure}
\caption{Multi-modal data integration on Archimedes Palimpsest. Top row: full experiment; bottom row: zoom. Parameters: $T=1\mathrm{e}5$, $\tau=1\mathrm{e}3$.}
\label{fig:archimede}
\end{figure}

The osmosis model can also be tailored to fuse image information with data coming from other image sources, see, e.g., in \cite{CalatroniUnveiling} an application to art restoration. In Figure \ref{fig:archimede}, we report the digital rendering of a CH image integrating the visible image with multi-modal pseudo colour data. 

Namely, having set $v_1$ to be the image in Figure \ref{fig: arch input} and $v_2$ to be the image in Figure \ref{fig: arch written} revealing the text underlying the Figure \ref{fig: arch input}, we post-process $v_2$ with a local Otsu thresholding method to segment the calligraphic data, while averaging the remaining information. Then we decompose the image domain $\Omega$ into $\Omega=\Omega_1\cup\Omega_2\cup\Omega_3$ where $\Omega_3$ is a tiny overlapping frame between an exterior rectangular frame $\Omega_1$ (the parchment) and an interior rectangular $\Omega_2$ (the text bounding box). Thus, we define the driving vector field $\bm{d}:=\bm{\nabla} \ln v$ with $v= \chi_{\Omega_1} v_1+ \chi_{\Omega_2} v_2+ \chi_{\Omega_3} (v_1+v_2)/2$ and run the osmosis model to get the composite result in Figure \ref{fig: arch written} integrating the pseudo-colour data visible on the input image, thus unveiling the hidden Archimedes' writing.

\section{Conclusions}\label{sec: conclusion}
In this paper we consider a variety of problems arising in CH imaging where the image osmosis model \eqref{eq:osmosis}  \cite{weickert,vogel,Calatroni2017} is used to correct and integrate multi-modal image acquisitions in different bands (from optical to thermal) in both reflectance and emission modes. To compute the numerical solution of the model efficiently for large image data, a dimensional reduction splitting approach is used. 

A free release of the MATLAB code used to compute the results reported in this work will be made available in the next future.

\subsubsection*{Acknowledgements.}  
SP acknowledges UK EPSRC grant EP/L016516/1. 
LC acknowledges the support of the Fondation Math\'ematique Jacques Hadamard (FMJH).
The diagnostics were supported by Dr.\ Francesca Tasso (\emph{Soprintendenza} of Castello Sforzesco) and Vittorio Barra (University of Verona).

\subsubsection*{Data Statement} 
The ``Monocromo'' data are sensitive: restricted access is subjected to the approval of ``Soprintendenza Castello, Musei Archeologici e Musei Storici'', Milan.
The Archimede palimpsest is released under CC-BY 3.0 license, \url{http://openn.library.upenn.edu/Data/0014/ArchimedesPalimpsest/}.

\end{document}